\documentclass[10pt]{article}
\usepackage[normalem]{ulem}
\usepackage{float,eurosym,cancel}
\usepackage{xspace,dsfont}
\usepackage[latin1]{inputenc}
\usepackage[dvipsnames]{xcolor}
\usepackage{amssymb}
\usepackage{amsmath}
\usepackage{amsthm}
\usepackage{enumerate}
\usepackage{graphicx}
\usepackage{caption}
\usepackage{subcaption}
\usepackage{multirow}
\usepackage{chicago}

\usepackage[textwidth = 155mm]{geometry}
\usepackage{tabularx}

\newtheorem{theorem}{Theorem}[section]
\newtheorem{lem}[theorem]{Lemma}

\newtheorem{cor}[theorem]{Corollary}
\newtheorem{conj}[theorem]{Conjecture}

\newtheorem{rem}[theorem]{Remark}

\newtheorem{ex}[theorem]{Example}
\newtheorem{defi}[theorem]{Definition}
\newtheorem{hyp}[theorem]{Assumption}

\newcommand{\bt}{\begin{theorem}}\newcommand{\et}{\end{theorem}}
\newcommand{\bl}{\begin{lem}}
\newcommand{\el}{\end{lem}}
\newcommand{\bcor}{\begin{cor}}\newcommand{\ecor}{\end{cor}}
\newcommand{\bconj}{\begin{conj}}\newcommand{\econj}{\end{conj}}

\newcommand{\bd}{\begin{defi} \rm }\newcommand{\ed}{\end{defi} }
\newcommand{\brem }{\begin{rem} \rm }\newcommand{\erem }{\end{rem}}
\newcommand{\bcom}{\begin{rem} \rm }\newcommand{\ecom }{\end{rem}}
\newcommand{\brems }{\begin{rem} \rm }\newcommand{\erems }{\end{rem}}
\newcommand{\bex}{\begin{ex} \rm }\newcommand{\eex}{\end{ex}}
\newcommand{\bhyp}{\begin{hyp} \rm }\newcommand{\ehyp}{\end{hyp}}

\def\bea{\begin{eqnarray}}
\def\eea{\end{eqnarray}}

\def\beqa{\begin{eqnarray}}
\def\eeqa{\end{eqnarray}}
\def\bal{\begin{aligned}}
\def\eal{\end{aligned}}

\def\Theta{{\mbox{CVA}}}

\begin{document}

\title{A general asymptotic formula for distinct partitions}
\author{Vivien Brunel $ ^{1,2}$
\\\\
$ ^1 $ L\'eonard de Vinci P\^ole Universitaire, Finance Lab, France \vspace{4pt}\\
$ ^2 $ Soci\'et\'e G\'en\'erale \vspace{4pt}
}

\date{\today}
\maketitle
\vspace{0cm}

\begin{abstract}
Many asymptotic formulas exist for unrestricted integer partitions as well as for distinct partitions of integers into a finite number of parts. Szekeres and Canfield have derived an asymptotic formula for the number of partitions that is valid for any value of the number of parts. We obtain general asymptotic formulas for distinct partitions that are valid in a wider range of parameters than the existing asymptotic formulas, and we recover the known asymptotic results as special cases of our general formulas. 

\vspace{0.2cm}

\end{abstract}
\newpage

\section*{Introduction}

Let $q(E,N)$ denote the number of partitions of an integer $E$ into $N$ distinct integer parts (these partitions are called distinct partitions),  and $P(E,N)$ the number of unrestricted integer partitions of $E$ into $N$ parts or less. It is well known in the theory of partitions that $q(E,N) = P(E-N(N+1)/, N)$. The asymptotic behaviour of $q(E,N)$ is often derived from the asymptotic behaviour of $P(E,N)$ and from the preceding relationship between the two.

The partition function of unrestricted partitions is known from the theory of partitions \citeN{Andrews1976}. A general asymptotic formula for $P(E,N)$ which is valid for all values of $N$ has been derived in \citeN{Szekeres1953}. \citeN{Canfield1997} has obtained it in a simpler way from the recursive formula satisfied by the function $P(E,N)$; however, his derivation requires to know the asymptotic dependence of the function $P(\cdot,N)$  to be derived.

Here, we solve the same problem for the number of distinct partitions. We obtain an explicit asymptotic formula for the function $q(E,N)$ when both $E$ and $N$ are large, which is similar to the one obtained by Canfield for unrestricted partitions. To our knowledge, this result is new. Additionally, we use the same approach to compute the asymptotic number $q(E,N,B)$ of distinct partitions with $N$ terms each one being smaller to to a given value $B$, in the limit when these three variables go to infinity but with constant values for $N^2/E = u$ and $N/B=p$. The asymptotic formula obtained by \citeN{Szekeres1987} corresponds to the limit $N^2/E \rightarrow \ln 2 / c$. Our formula is valid in a wider range of parameters since $u$ is arbitrary; we recover Szekeres' formula as a special case.

The plan of the paper is as follows: we first detail our methodology.Second, we derive Szekeres-Canfield's formula for unrestricted partitions with a saddle-point approximation which, in essence, is similar to Szekeres' derivation. Third, we derive the general asymptotic formula for the number of unequal partitions $q(E,N)$. Finally, we extend the formula to cases where, in addition, the summands constrained to be smaller or equal to a given integer $B$.  

\section{Methodology and notations}

Throughout this paper we follow the same methodology. We illustrate it with the function $P(E,N)$. We introduce the partition function related to the function $P(E,N)$:
\begin{equation}
Z_E(x,z) = \sum_{E=1}^{\infty} \sum_{N=1}^{\infty} z^{N} \, x^{E} P(E,N)
\end{equation}

So, by inverting this relationship, we have:
\begin{equation}{\label{curviligne}}
P(E,N) = \frac{1}{(2i \pi)^2} \oint \oint dx \, dz  \frac{ Z_E(x,z)}{z^{N+1} \, x^{E+1}}
\end{equation}
We set $z=e^{-\alpha}$ and $x=e^{-\beta}$.
The partition function for unrestricted partitions is:
\begin{equation}
Z_E(\alpha, \beta) = \prod_{k=1}^{\infty} \left(  1-e^{-\alpha}e^{-\beta k} \right )^{-1}
\end{equation}
 The partition function of $P(E,N)$ is well known \cite{Andrews1976}:
\begin{equation}{\label{PartitionFunction1}}
Z_E(\alpha, \beta) = \prod_{k=1}^{\infty} \left(  1-e^{-\alpha} e^{-\beta k} \right )^{-1}
\end{equation}

The integrals of eq.(\ref{curviligne}) are taken over closed contours round the origin of the complex x and z plane. We apply the saddle-point approximation method to compute eq.(\ref{curviligne}). 
So, by inverting this relationship, we have:
\begin{equation}
P(E,N) \sim \frac{1}{2 \pi \sqrt{ D}} \exp \left [ \alpha^* N + \beta^* E + \ln Z_E(\alpha^*, \beta^*)
  \right ] 
\end{equation}
where $\alpha^*$ and $\beta^*$ are the coordinates of the saddle-point that maximizes the entropy function $S_E(\alpha, \beta) = \alpha N + \beta E + \ln Z_E(\alpha, \beta)$ and:
\begin{equation}
D_E = \left ( \frac{\partial^2 S_E}{\partial \beta^2}  \right) . \left ( \frac{\partial^2 S_E}{\partial \alpha^2}  \right) - \left ( \frac{\partial^2 S_E}{\partial \beta \partial \alpha}  \right)^2
\end{equation}
The second order derivatives are evaluated at the saddle-point $(\alpha^*, \beta^*)$. 

Similar formulas hold for the number of distinct partitions $q(E,N)$  and $q(E,N,B)$. The partition function of $q(E,N)$ is:
\begin{equation}
Z_D(\alpha, \beta) = \prod_{k=1}^{\infty} \left(  1-e^{-\alpha}e^{-\beta k} \right )^{-1}
\end{equation}
The associated entropy function is $S_D(\alpha, \beta) = \alpha N + \beta E + \ln Z_D(\alpha, \beta)$. The restriction of the size of  parts to be lower than $B$ modifies the partition function of $q(E,N,B)$ compared to the one of $q(E,N)$:
\begin{equation}
Z_R(\alpha, \beta) = \prod_{k=1}^{B} \left(  1-e^{-\alpha}e^{-\beta k} \right )^{-1}
\end{equation}
The entropy function $S_R(\alpha, \beta)$ is changed accordingly.

Finally, following \citeN{Szekeres1987}, we introduce the following notations:
\begin{equation}
\begin{array}{l}
c^2 = \frac{\pi^2}{12}\\
\gamma = 1-\frac{(\ln 2)^2}{c^2} \\
\sigma = N-\sqrt{E}\frac{\ln 2}{c} 
\end{array}
\end{equation}

\section{Canfield's formula from the saddle-point approximation}

We recover Canfield's formula for $P(E,N)$ in an easy way with the saddle-point approximation.
The entropy function writes for equal partitions:
\begin{eqnarray}
S_E (\alpha, \beta) &=&  \alpha N + \beta E  - \sum_{k=1}^{\infty} \ln(1- e^{-\alpha} e^{-\beta k} ) \\
&=&  \alpha N +\beta E - \sum_{k=0}^{\infty} \ln(1- e^{-\alpha} e^{-\beta k} ) + \ln(1-e^{-\alpha})
\end{eqnarray}

By applying the Euler-McLaurin formula to the sum term, we get:
\begin{equation}
S_E (\alpha, \beta) \sim  \alpha N + \beta E - \frac{1}{\beta} \int_{0}^{\infty} \ln(1 -e^{-\alpha} e^{-x} ) dx + \frac{1}{2} \ln(1-e^{-\alpha})
\end{equation}

The saddle-point equations write:
\begin{equation}{\label{saddleB2}}
\frac{\partial S_E}{\partial \alpha} = 0 \sim N - \frac{1}{\beta} \int_{0}^{\infty} \frac{dx}{e^{\alpha+x}-1} 
\end{equation}
and 
\begin{equation}{\label{saddleB1}}
\frac{\partial S_E}{\partial \beta}  = 0 \sim E + \frac{1}{\beta^2} \int_{0}^{\infty} \ln (1-e^{-\alpha}e^{-x}) dx
\end{equation}

We set $v=\beta^* N$ and $u=N/\sqrt{E}$
After some analytics, we get the value of $\alpha$ at the saddle-point:
\begin{equation}{\label{saddle2B}}
e^{-\alpha^*} = 1-e^{-\beta^* N} =  1 - e^{-v}
\end{equation}
The integral term in eq.(\ref{saddleB1}) is equivalent to:
\begin{equation}
\int_{0}^{\infty} \ln(1- e^{-\alpha^*} e^{-x} ) dx  = \int_0^{\beta^* N} \frac{t\, dt}{1-e^{t}}
\end{equation}
Then, eq.(\ref{saddleB1}) is equivalent to:
\begin{equation}
\frac{v^2}{u^2}  =  \int_0^{v} \frac{t\, dt}{e^{t}-1}
\end{equation}
The second order derivatives are:
\begin{equation}
\begin{array}{l}
\frac{\partial^2 S_E}{\partial \beta^2} \sim  - \frac{2}{\beta^{*3}} \int_{0}^{\infty} \ln (1-e^{-\alpha^*}e^{-x}) dx =  \frac{2}{\beta^{*3}} \frac{v^2}{u^2} \\
\frac{\partial^2 S_E}{\partial \alpha \partial \beta} \sim  \frac{1}{\beta^{*2}} \frac{d}{d\alpha} \int_{0}^{\infty} \ln (1-e^{-\alpha^*}e^{-x}) dx = \frac{v}{\beta^{*2}} \\
\frac{\partial^2 S_E}{\partial \alpha^2} \sim \frac{e^{v}-1}{\beta^*}
\end{array}
\end{equation}
At the saddle-point, the  function $S_E (\alpha^*, \beta^*)$ is equal to:
\begin{eqnarray}
S_E (\alpha^*, \beta^*) & =& 2 \beta^* E + \alpha^* N + v/2 \\
&=& \left [  2\frac{v}{u} - u \ln (1-e^{-v}) \right ] \sqrt{E} + \frac{v}{2}
\end{eqnarray}
The denominator of the saddle-point formula is equal to:
\begin{eqnarray}
2 \pi \sqrt{D_E} &=& \left ( \left( \frac{2}{\beta^{*3}} \frac{v^2}{u^2}  \right) \cdot \left ( \frac{e^{v}-1}{\beta^*} \right)  -  \frac{v^2}{\beta^{*4}} \right)^{1/2} \\
 &=& 2^{3/2} \pi E \frac{u}{v} e^{v/2} \left ( (1-e^{-v}) - \frac{1}{2} u^2 e^{-v}  \right)^{1/2}
\end{eqnarray}

Finally, we get the number of  partitions of $E$ into  $N$ summands or less:
\begin{equation}{\label{Canfield}}
P(E,N) = \frac{f(u)}{E} \, e^{ \sqrt{E}g(u)}
\end{equation}
where:
\begin{equation}
g(u) =  2\frac{v}{u} - u \ln (1-e^{-v}) 
\end{equation}
and
\begin{equation}
f(u) =  \frac{v}{2^{3/2} \pi u}  \left ( (1-e^{-v}) - \frac{1}{2} u^2 e^{-v}  \right)^{-1/2}
\end{equation}
This is exactly Canfield's formula \cite{Canfield1997}.

\section{Asymptotic formula for unrestricted partitions}

\subsection{General formula for $q(E,N)$}

We use the same methodology to obtain a similar formula for unequal partitions. The entropy formula for unequal partitions writes:
\begin{eqnarray}{\label{entropyfunction}}
S_D (\alpha, \beta) &=& \beta E + \alpha N + \sum_{k=1}^{\infty} \ln(1+ e^{-\alpha} e^{-\beta k} ) \\
&=& \beta E + \alpha N + \sum_{k=0}^{\infty} \ln(1+ e^{-\alpha} e^{-\beta k} ) - \ln(1+e^{-\alpha})
\end{eqnarray}

By applying the Euler-McLaurin formula to the sum term, we get:
\begin{equation}
S_D (\alpha, \beta) \sim \beta E + \alpha N + \frac{1}{\beta} \int_{0}^{\infty} \ln(1+ e^{-\alpha} e^{-x} ) dx - \frac{1}{2} \ln(1+e^{-\alpha})
\end{equation}

The saddle-point equations write:
\begin{equation}{\label{saddle2f}}
\frac{\partial S_D}{\partial \alpha} = 0 \sim N - \frac{1}{\beta} \int_{0}^{\infty} \frac{dx}{1+e^{\alpha + x}} 
\end{equation}
and
\begin{equation}{\label{saddle1f}}
\frac{\partial S_D}{\partial \beta}  = 0 \sim E - \frac{1}{\beta^2} \int_{0}^{\infty} \ln (1+e^{-\alpha}e^{-x}) dx
\end{equation}

The coordinates $(\alpha^*, \beta^*)$ of the saddle-point for distinct partitions are different from the ones for unrestricted partitions as seen in the previous section.
We set $v=\beta^* N$ and $u=N/\sqrt{E}$
After some analytics, we get the value of $\alpha^*$ at the saddle-point:
\begin{equation}{\label{saddle2f2}}
e^{-\alpha^*} = e^{\beta^* N} - 1 = e^v - 1
\end{equation}
 The integral term in eq.(\ref{saddle1f}) is equivalent to:
\begin{equation}
\int_{0}^{\infty} \ln(1+ e^{-\alpha^*} e^{-x} ) dx  = \int_0^{\beta^* N} \frac{t\, dt}{1-e^{-t}}
\end{equation}
Then, eq.(\ref{saddle1f}) is equivalent to:
\begin{equation}{\label{saddle1f2}}
\frac{v^2}{u^2}  = \int_0^{v} \frac{t\, dt}{1-e^{-t}}
\end{equation}

The second order derivatives are:
\begin{equation}
\begin{array}{l}
\frac{\partial^2 S_D}{\partial \beta^2} \sim  \frac{2}{\beta^{*3}} \int_{0}^{\infty} \ln (1+e^{-\alpha^*}e^{-x}) dx = \frac{2}{\beta^{*3}} \frac{v^2}{u^2} \\
\frac{\partial^2 S_D}{\partial \alpha \partial \beta} \sim - \frac{1}{\beta^{*2}} \frac{d}{d\alpha} \int_{0}^{\infty} \ln (1+e^{-\alpha^*}e^{-x}) dx = \frac{v}{\beta^{*2}} \\
\frac{\partial^2 S_D}{\partial \alpha^2} \sim \frac{1-e^{-v}}{\beta^*}
\end{array}
\end{equation}

At the Saddle-point, the entropy function is equal to:
\begin{eqnarray}
S_D (\alpha^*, \beta^*) & =& 2 \beta^* E + \alpha^* N - v/2 \\
&=& \left [  2\frac{v}{u} - u \ln (e^v-1) \right ] \sqrt{E} - \frac{v}{2}
\end{eqnarray}

The denominator of the saddle-point formula is equal to:
\begin{eqnarray}
2 \pi \sqrt{D_D} &=& \left ( \left( \frac{2}{\beta^{*3}} \frac{v^2}{u^2}  \right) \cdot \left ( \frac{1-e^{-v}}{\beta^*} \right)  -  \frac{v^2}{\beta^{*4}} \right)^{1/2} \\
 &=& 2^{3/2} \pi E \frac{u}{v} \left ( (1-e^{-v}) - \frac{1}{2} u^2  \right)^{1/2}
\end{eqnarray}

Finally, we get the number of unequal partitions of $E$ into $N$ summands:
\begin{eqnarray}
q(E,N) = \frac{\bar{f}(u)}{E} \, e^{ \sqrt{E}\bar{g}(u)}
\end{eqnarray}
where:
\begin{equation}{\label{denomF}}
\begin{array}{l}
\bar{g}(u) =  2\frac{v}{u} - u \ln (e^v-1) \\
\bar{f}(u) =  \frac{v}{2^{3/2} \pi  u  }   \left ( (e^{v}-1) - \frac{1}{2} u^2 e^{v} \right)^{-1/2} \\
\frac{v^2}{u^2}  = \int_0^{v} \frac{t\, dt}{1-e^{-t}}
\end{array}
\end{equation}

\subsection{Limit $N << \sqrt{E}$}

Let's now focus on some limiting cases.
In the limit $N << \sqrt{E}$ (i.e. $u<<1$), we have $\alpha^* \rightarrow \infty$, and then, from eq.(\ref{saddle2f2}), $v=\beta^*  N \rightarrow 0$. From eq.(\ref{saddle1f}):
\begin{equation}
\beta^{*2} E = \int_0^{\infty} \ln (1+e^{-\alpha^*}e^{-x}) dx \sim e^{-\alpha^*}
\end{equation}
Then, $\beta^* = \sqrt{e^{-\alpha^*}/E}\sim \sqrt{\beta^* N/E}$ because $\beta^* N <<1$. It follows that $\beta^* \sim N/E$. We obtain $u \sim \sqrt{\beta^* N}$ and $v \sim \beta^* N$. The function $\bar{g}(u)$ is equal to:
\begin{equation}
\bar{g}(u) = 2 \frac{N}{\sqrt{E}} - \frac{N}{\sqrt{E}} \ln \left ( \frac{N^2}{E}  \right)
\end{equation}
Similarly, from eq.(\ref{denomF}), we get for $\bar{f}(u)$:
\begin{equation}
\bar{f}(u) =  \frac{1}{2^{3/2} \pi  \frac{N}{\sqrt{E}}  e^{N^2/2E} \left (  1-  e^{-N^2/E} - \frac{N^2}{2E} \right)^{1/2}}
\end{equation}
We obtain the known result:
\begin{equation}{\label{LDPlimit}}
q(E,N) \sim \frac{e^{2N} E^{N-1}}{2 \pi N^{2N}}
\end{equation}
We can take the same limit  $N << \sqrt{E}$ in Canfield's formula eq.(\ref{Canfield}) for unrestricted partitions and we will get the same result as in eq.(\ref{LDPlimit}). Both asymptotics converge to the same limiting distribution which corresponds to the Maxwell-Boltzman distribution as explained for instance in \cite{Debnath1987}.

\subsection{Limit $N \rightarrow \sqrt{2E}$}

We know from eq.(\ref{saddle2f2}), that $N \leq \sqrt{2E}$. We can compute the asymptotics of $q(E,N)$ in the limit $N \rightarrow \sqrt{2E}$. This limit correspond to $\alpha^*\rightarrow -\infty$ and $\beta^* N \rightarrow \infty$.
It is straightforward to see that $\beta^{*2} (E-N^2/2) \sim 2c^2$:
\begin{eqnarray*}
\beta^{*2} (E-N^2/2)  &=& \int_0^{\infty} \ln (1+e^{\beta^*  N}e^{-x}) dx - \int_0^{\beta^*  N} x \, dx \\
&=&  \int_0^{\infty} \ln (1+e^{-x}) dx + \int_{-\beta^* N}^{0} \ln (1+e^{-x}) dx- \int_0^{\beta^* N} x \, dx \\
&=& c^2 + \int_{0}^{\beta^* N} \ln (1+e^{x}) -x  dx \\
&=& c^2 + \int_{0}^{\beta^* N} \ln (1+e^{-x}) dx \\
&\sim& 2c^2
\end{eqnarray*}

\subsection{Limit $\alpha^* \rightarrow 0$}

Finally, there is another interesting limit when $\alpha^* \rightarrow 0$. This limit corresponds to the total number of unequal partitions, whatever the value of $N$.

In the limit $\alpha^* \rightarrow 0$, we have $v\sim \ln 2$. We set $v = \ln2 + \epsilon$. At first order in $\epsilon$, we obtain from eq.(\ref{saddle1f2}):
\begin{equation}{\label{approximations}}
\begin{array}{l}
\frac{v}{u}= c + \frac{\ln 2}{c} \epsilon \\
u \ln (e^v - 1) = \frac{2 \ln2}{c} \epsilon + \frac{2 \gamma}{c} \epsilon^2
\end{array}
\end{equation}
Then, we obtain $g(u) \sim 2c - 2c\gamma \epsilon^2/c$
In addition, we can rewrite $\epsilon$ as a function of $N$ and $E$. Indeed:
\begin{eqnarray}
\epsilon &=& v- \ln2 = u\frac{v}{u} - \ln2 \\
& =& \frac{N}{\sqrt{E}} \left(  c+ \frac{\ln2}{c} \epsilon \right ) - \ln2
\end{eqnarray}
We can solve this first order equation to get $\epsilon \sim \frac{1}{\gamma} \left (  c\frac{N}{\sqrt{E}} - \ln 2\right )$. Finally, we obtain
\begin{equation}
\epsilon^2 \sim \frac{c^2 \sigma^2}{\gamma^2 E}
\end{equation}
From the approximations of eq.(\ref{approximations}), we get $\bar{f}(u) \sim 4\sqrt{6 \gamma}$. In the limit $\alpha^* \rightarrow 0$, the number of unequal partitions with N terms is:
\begin{equation}
q(E,N) \sim \frac{\exp \left (  2c\sqrt{E} - \frac{2c\sigma^2}{\gamma \sqrt{E}} \right )}{4\sqrt{6 \gamma}}
\end{equation}
This is exactly the result of \citeN{Erdos1941}. \citeN{Szekeres1987} explains that the total number of unequal partitions is obtained by summing over all the values of $N$ (or equivalently over all values of $\sigma$), and we obtain the well known asymptotic formula $q(E) \sim \exp (  2c\sqrt{E}) / \left(  4. 3^{1/4}. E^{3/4} \right)$.

\section{Formula when the summands are smaller than $B$}

In the case where the summands are restricted to be lower than a given number $B$, the entropy function of eq.(\ref{entropyfunction}) is slightly modified because the sum term stops at $B$ instead of infinity:
\begin{equation}{\label{entropy}}
S_R(\alpha, \beta) = \alpha N + \beta E +  \sum_{k=0}^{B} \ln(1+ e^{-\alpha} e^{-\beta k} ) - \ln(1+e^{-\alpha})  
\end{equation}

We introduce the partition function related to the function $q(E,N,B)$:
\begin{equation}
Z_R(x,z) = \sum_{E=1}^{\infty} \sum_{N=1}^{\infty} z^{N} \, x^{E} q(E,N,B)
\end{equation}

We set $z=e^{-\alpha}$ and $x=e^{-\beta}$.
 The partition function for unequal partitions is well known:
\begin{equation}{\label{PartitionFunction1}}
Z_R(\alpha, \beta) = \prod_{k=1}^{B} \left(  1+e^{-\alpha} e^{-\beta k} \right )
\end{equation}

By inverting this relationship, we have:
\begin{equation}
q(E,N,B) \sim \frac{1}{2 \pi \sqrt{ D_R}} \exp \left [ \alpha^* N + \beta^* E + \ln Z_R(\alpha^*, \beta^*)
  \right ] 
\end{equation}
where $\alpha^*$ and $\beta^*$ are the coordinates of the saddle-point that maximizes the function $S_R(\alpha, \beta) = \alpha N + \beta E + \ln Z_R(\alpha, \beta)$ and:
\begin{equation}
D_R = \left ( \frac{\partial^2 S}{\partial \beta^2}  \right) . \left ( \frac{\partial^2 S}{\partial \alpha^2}  \right) - \left ( \frac{\partial^2 S}{\partial \beta \partial \alpha}  \right)^2
\end{equation}
The second order derivatives are taken at the saddle-point $(\alpha^*, \beta^*)$.

By applying the Euler-McLaurin formula to the sum term of eq.(\ref{entropy}), we get:
\begin{equation}
S_R(\alpha, \beta) = \beta E + \alpha N + \frac{1}{\beta} \int_{0}^{\beta B} \ln(1+ e^{-\alpha} e^{-x} ) dx - \frac{1}{2} \ln(1+e^{-\alpha}) + \frac{1}{2} \ln(1+e^{-\alpha -\beta B}) 
\end{equation}

The saddle-point equations write:
\begin{equation}{\label{saddle2F}}
\frac{\partial S}{\partial \alpha} = 0 = N - \frac{1}{\beta} \int_{0}^{\beta B} \frac{dx}{1+e^{\alpha + x}} 
\end{equation}
and
\begin{equation}{\label{saddle1F}}
\frac{\partial S}{\partial \beta}  = 0 = E - \frac{1}{\beta^2} \int_{0}^{\beta B} \ln (1+e^{-\alpha}e^{-x}) dx + \frac{B}{\beta}  \ln (1+e^{-\alpha- \beta B}) 
\end{equation}

We set $p=N/B$, $v=\beta^* N$, $e^{-\alpha^*} = e^{v(p)}-1$ and $u=N/\sqrt{E}$.
After some analytics, we get the value of $\alpha^*$ at the saddle-point from eq.(\ref{saddle2F}):
\begin{equation}{\label{saddle2F2}}
e^{-\alpha^*} = e^{v(p)}-1 = \frac{e^{v} - 1}{1-e^{v(1-1/p)}}
\end{equation}
Additionally, eq.(\ref{saddle2F}) leads to $e^{-\alpha- \beta B} \sim v(p)-v$.
 The integral term in eq.(\ref{saddle1F}) is equal to:
\begin{equation}
\int_{0}^{\beta B} \ln(1+ e^{-\alpha} e^{-x} ) dx  = \int_0^{\ln (1+ e^{-\alpha})} \frac{t\, dt}{1-e^{-t}} - e^{-\alpha - \beta B}
\end{equation}
We see that the finite size correction coming from the upper bound of the integral is much smaller than the last term in eq.(\ref{saddle1F}), but it cannot be neglected as we shall see later. Then, eq.(\ref{saddle1F}) is equivalent to:
\begin{equation}{\label{saddlepoint}}
\frac{v^2}{u^2}  = \int_0^{v(p)} \frac{t\, dt}{1-e^{-t}} + \left( \frac{v}{p}-1 \right) \left ( v(p)-v \right )
\end{equation}

The second order derivatives are, when $B \rightarrow \infty$:
\begin{equation}
\begin{array}{l}
\frac{\partial^2 S_R}{\partial \beta} \sim  \frac{2}{\beta^{*3}} \int_{0}^{\infty} \ln (1+e^{-\alpha^*}e^{-x}) dx = \frac{2}{\beta^{*3}} \frac{v^2}{u^2} \\
\frac{\partial^2 S_R}{\partial \alpha \partial \beta} \sim - \frac{1}{\beta^{*2}} \frac{d}{d\alpha} \int_{0}^{\infty} \ln (1+e^{-\alpha^*}e^{-x}) dx = \frac{v}{\beta^{*2}} \\
\frac{\partial^2 S_R}{\partial \alpha^2} \sim \frac{1-e^{-v(p)}}{\beta^*}
\end{array}
\end{equation}

At the saddle-point, the function $S(\alpha^*, \beta^*)$ is equal to:
\begin{eqnarray}
S_R (\alpha^*, \beta^*) & =& 2 \beta^* E + \alpha^* N  + \frac{1}{\beta^*} \frac{v}{p} \left ( v(p)-v \right ) - v(p)/2\\
&=& \left [  2\frac{v}{u} - u \ln (e^{v(p)}-1) +\frac{u}{p} \left ( v(p)-v \right ) \right ] \sqrt{E} - \frac{v(p)}{2} 
\end{eqnarray}

The denominator of the saddle-point formula is equal to:
\begin{eqnarray}
2 \pi \sqrt{D_R} &=& \left ( \left( \frac{2}{\beta^{*3}} \frac{v^2}{u^2}  \right) \cdot \left ( \frac{1-e^{-v(p)}}{\beta^*} \right)  -  \frac{v^2}{\beta^{*4}} \right)^{1/2} \\
 &=& 2^{3/2} \pi E \frac{u}{v} \left ( (1-e^{-v(p)}) - \frac{1}{2} u^2  \right)^{1/2}
\end{eqnarray}

Finally, we get the number of unequal partitions of $E$ into $N$ summands, each summand being not larger than $B$:
\begin{eqnarray}{\label{eqgenerale1}}
q(E,N,B) = \frac{\tilde{f}(u)}{E} \, e^{ \sqrt{E}\tilde{g}(u)}
\end{eqnarray}
where:
\begin{equation}{\label{eqgenerale2}}
\begin{array}{l}
\tilde{g}(u) =  2\frac{v}{u} - u \ln (e^{v(p)}-1) + \frac{u}{p} \left ( v(p)-v \right ) \\
\tilde{f}(u) =  \frac{v}{2^{3/2} \pi  u  }   \left ( (e^{v(p)}-1) - \frac{1}{2} u^2 e^{v(p)} \right)^{-1/2}
\end{array}
\end{equation}

In the limit $\alpha \rightarrow 0$, we have $v(p), \, \,  v \sim \ln 2$. We set $v(p) =  \ln2 + \epsilon + (v(p)-v)$. At first order in $\epsilon$ and $v(p)-v$, we obtain from eq.(\ref{saddlepoint}):
\begin{equation}{\label{approximations2}}
\begin{array}{l}
\frac{v}{u} = c + \frac{\ln 2}{c} \epsilon + \frac{1}{2c} \left(  2 \ln 2- \frac{v}{p} -1 \right)  (v(p)-v)\\
u \ln (e^{v(p)} - 1) = \frac{2 \ln2}{c} \epsilon + 2c \gamma \epsilon^2 + \frac{2 \ln2}{c} (v(p)-v)
\end{array}
\end{equation}
Then, we obtain $\tilde{g}(u) \sim 2c - 2 \gamma \epsilon^2/c - \frac{1}{c} (v(p)-v)$.
Compared to the case of $q(E,N)$, there is an additional term in the definition of the function $\tilde{g}(u)$. This leads directly to Szekeres' asymptotic formula as a special case of eq.(\ref{eqgenerale1} - \ref{eqgenerale2}):
\begin{equation}
q(E,N,B) \sim \frac{\exp \left (  2c\sqrt{E} - \frac{2c\sigma^2}{\gamma \sqrt{E} } - \frac{\sqrt{E}}{c} e^{cB/\sqrt{E}} \right )}{4\sqrt{6 \gamma}E}
\end{equation}

\vspace{1cm}
\noindent Acknowledgement: I am grateful to Denis Bernard for the fruitful discussions we had together.

\end{document}